\documentclass[a4paper,11pt]{article}
\usepackage[latin1]{inputenc}
\usepackage[active]{srcltx}
\usepackage{amsfonts}
\usepackage{amsmath}
\usepackage{amssymb}
\usepackage{amsthm}

\usepackage{chapterbib}
\usepackage{bibunits}
\usepackage[sectionbib]{natbib}
\usepackage{bm}
\usepackage{subfigure}

\usepackage[active]{srcltx}
\usepackage{ifpdf}
\usepackage{setspace}
\usepackage{lscape}
\usepackage{bm}

\usepackage{mathrsfs}
\usepackage{relsize}

\usepackage{color}
\definecolor{darkred}{rgb}{0.5,0,0}
\definecolor{darkgreen}{rgb}{0,0.5,0}
\definecolor{darkblue}{rgb}{0,0,0.5}

\ifpdf

\newcommand{\drivy}{pdftex}
\usepackage[pdftex=true,hyperfootnotes=true]{hyperref}
\else

\newcommand{\drivy}{dvips}
\usepackage[hyperfootnotes=true]{hyperref}
\fi
\usepackage[\drivy]{graphicx}
\usepackage{rotating}

\usepackage{setspace}
\usepackage[\drivy]{graphicx}
\usepackage[footnotesize,nooneline,ignoreLTcapwidth,bf]{caption2}

\theoremstyle{definition}

\newtheorem{remark}{Remark}

\date{}

\newcommand{\vnu}{\bm \nu}

\makeatletter
\newcommand{\rd}{\@ifnextchar^{\DIfF}{\DIfF^{}}}
\def\DIfF^#1{%
   \mathop{\mathrm{\mathstrut d}}%
   \nolimits^{#1}\gobblespace}
\def\gobblespace{\futurelet\diffarg\opspace}
\def\opspace{%
   \let\DiffSpace\!%
   \ifx\diffarg(%
   \let\DiffSpace\relax
   \else
   \ifx\diffarg[%
   \let\DiffSpace\relax
   \else
   \ifx\diffarg\{%
   \let\DiffSpace\relax
   \fi\fi\fi\DiffSpace}

\usepackage{pdfpages}

\title{A new strategy for Robbins' problem of optimal stopping\thanks{
The authors thank F. Thomas Bruss,  Alexander Gnedin, Klaus P\"otzelberger and Ester Samuel-Cahn for interesting discussions and comments.
We thank the referee and especially the associate editor for their very helpful comments and suggestions. The suggestions of the associate editor improved the exposition a lot.
}}
\author{Martin Meier  and Leopold S\"ogner%
{\thanks{Martin Meier, 
(meier@ihs.ac.at), Leopold S\"ogner, (soegner@ihs.ac.at), Department of Economics and Finance, 
Institute for Advanced Studies, 
Josefst\"adter Stra\ss{}e 39,
1080 Vienna, Austria. 
Martin Meier has a further affiliation with the Vienna Graduate School of Economics (VGSE).
Leopold S\"ogner has a further affiliation with the Vienna Graduate School of Finance (VGSF).
}
}} 

\date{\today}

\begin{document}
\maketitle

\begin{abstract}
In this article we study the expected rank problem under full information. 
Our approach uses the planar Poisson approach from \citet{Gnedin2007} to derive the expected rank of a stopping rule that is one of the simplest non-trivial examples combining rank dependent rules with threshold rules.  This rule attains an expected rank lower than the best upper bounds obtained in the literature so far, in particular we obtain an expected rank of $2.32614$. 

\vskip1mm
\noindent\textit{Keywords:} Optimal stopping, Robbins' problem\\[-5pt] 
\vskip-2mm
\noindent\textit{2000 Mathematics Subject Classification:} Primary 60G40

\end{abstract}

\section{Introduction and Motivation}
\label{sect:intro}

Consider $iid$ distributed random variables $X_1,\dots,X_n$. The rank of $X_k$ is defined as $R_k = \sum_{m=1}^{n} \mathbf{1}_{(X_m \leq X_k)}$.  Robbins' problem is to find the optimal stopping rule that minimizes the expected rank of the chosen realization.  
 Although the optimal rule is unknown, some results on the properties of the optimal rule are available: \citet{Bruss96halfprophet}[Section~4.2] proved full history dependence of the optimal rule for 
the $n$-period problem. That is to say, the optimal decision to stop at stage $m<n$ depends on all past realizations
$x_1,\dots,x_{m-1}$ as well as the current realization $x_m$.   
In addition,  by using variational calculus \citet{AssafChan1996} obtained a lower bound of $1.85$ for the limit of the $n$-period problem.
Using truncated loss functions, \citet{BrussFerguson1993} derived a lower bound of $1.908$ by computational methods (with the help of Hardwick and Schork). 
A further question raised by \citet{BrussSwan2009}, which is associated with Robbins' problem, is whether the limit superior of $n$ times the expected value of $X_{\tau_n}$ arising from the optimal $n$-period rule $\tau_n$
is finite. This question was settled by \citet{GnedinIksanov2011}, who showed that it is indeed finite. 
For an overview on Robbins' problem and related stopping problems the reader is referred to \citet{Bruss2005}, \citet{Gnedin2007} and \citet{Swan2011}. 

{\em Full information rules} are stopping rules adapted to the filtration generated by all prior and current observations.
Interesting subsets of full information rules are {\em rank rules}, where the decision to stop depends on time $k$ and the relative rank $I_k = \sum_{m=1}^{k} \mathbf{1}_{(X_m \leq X_k)}$, 
 and  
{\em threshold rules}, where one stops with the first realization $x_m$ such that  $x_m \leq f_m$, where $f_m$ is some positive real number that depends only on $m$.
When minimizing the expected rank under full information for the $n\rightarrow \infty$ case, \cite{AssafChan1996}[Example~4.2] showed that a value of $2.3318$ can be obtained with a threshold rule of a simple form. This value was replicated in \cite{Gnedin2007} for the continuous time Poisson embedding.
A value of $7/3$ has been obtained in \cite{BrussFerguson1993} with another approximately optimal threshold rule. 
They estimated that the optimal threshold rule gives a value of $2.32659$. Additionally, \cite{AssafChan1996} showed that in the $n \rightarrow \infty$-limit the expected loss  for the optimal threshold rule has to be in the interval $(2.295,2.3267)$.

More recent literature such as \cite{Gnedin1996,Gnedin2004,Gnedin2007}, \cite{Bruss2001}, and \cite{BrussSwan2009} has shifted attention to Poisson embeddings of discrete time optimal stopping problems.  
These articles also demonstrate how the continuous time versions can be used to obtain upper bounds for the $n$-period discrete time problems. 
By using a continuous time Poisson embedding of Robbins' problem, \cite{Gnedin2007} showed that full history dependence also persists in the $n \rightarrow \infty$ limit, so that even in the limit a threshold rule cannot be optimal. However it is still not clear by how much the optimal rule is better than the optimal threshold rule.

This article combines a threshold rule with a rank dependent rule, which enables us to obtain an analytic solution with an expected rank smaller than $2.32659$. More specifically, up to a time $\alpha$ we use a rule like \cite{Gnedin2004} for the best choice problem, where one stops with the first observation that is below some threshold function $f_b(t)$ and has relative rank $1$, 
but with a different threshold parameter $b$. From $\alpha$ on, we apply a threshold rule, where -- given that the stopping criterion was not fulfilled for $t \leq \alpha$ -- the decision maker stops at the first observation below a function $f_c(t)$, again with a different threshold parameter $c$ as in \cite{Gnedin2007}[Section~3]. The framework of \cite{Gnedin2007} allows us to work ``directly in the limit". It would be much harder to use a similar family of rules like our rule in the $n$-period problems, and then compute the expected rank of this family of rules as $n \rightarrow \infty$. 

A similar rule for the discrete time problem has already been proposed by \cite{AssafChan1996}[Remark~6.2]. There, an approximately optimal threshold rule is combined with the requirement to stop with a relative rank of one within a time span $\gamma  n$, where $\gamma \in [0,1)$. 
An analytical investigation of this rule -- in the planar Poisson framework of \cite{Gnedin2004} -- was performed in \cite{Tamaki2004}, who obtains a value of $2.33044$, with $\alpha=0.42$ and $b=c=1.95$.  The main difference to our approach is that the same threshold function is applied during the whole time span, which simplifies the computations compared to ours'.



\section{A Simple Rank-Threshold Rule}
\label{sect:martin1}

We follow 
\cite{Gnedin2004,Gnedin2007} and consider a scatter of atoms $\mathscr{P}$ arising from a continuous time planar Poisson process
on the strip $[0,1] \times
\bar{\mathbb{R}}_{+}$, where $\bar{\mathbb{R}}_{+}$ stands for the interval $[0,\infty]$.
The intensity measure is the Lebesgue measure $dt dx$, 
which implies that the number of particles $N$ in some subset with Lebesgue measure $\vnu$ 
follows a Poisson distribution with density $\mathbb{P}(N=n)= e^{- \vnu}  \frac{\vnu^n }{n!}$. 
An atom $(T,X)$ consists of the arrival time $T$ 
and the value $X$. 
%
%
Ordering the atoms with respect to $X_j$ in ascending order yields 
the increasing sequence of points $\left( X_{1,1},X_{1,2},X_{1,3}, \dots \right)=: \mathbf{X}_1$ of a unit Poisson process. We denote an atom of the ordered sequence by $(T_r,X_{1,r})$. By the properties of the planar Poisson process the arrival times $T_r$ are uniform $iid$ on $[0,1]$, and $T_r$, $r=1,2,\dots$, and $\mathbf{X}_1$ are independent. 
In addition, a 
stochastic process $\left( \tilde{X}_t,\ t \in [0,1] \right)$ with values in $[0,\infty] $ can be constructed, by $\tilde{X}_t := X_r$, where $r$ is the minimal $r$ such that $t = T_r$, if such an $r$ exists, and else 
$\tilde{X}_t := \infty$.    
For more technical details on this planar Poisson process the reader is referred to \cite{Gnedin2004,Gnedin2007}. 

Let $(t,x)$ and $\mathcal{P}$ stand for realizations of $(T,X)$ and $\mathscr{P}$. 
For a generic $(t ,x) \in \mathcal{P}$ the absolute rank, $R(t ,x)$, is defined by the number of $\mathcal{P}$-points strictly south of $(t,x)$
plus $1$, while the relative rank, $I(t ,x)$, is the number of $\mathcal{P}$-points strictly south-west of $(t, x)$ plus $1$.

Let $(\tau,X)$ define a {\em stopping point}, where the event 
$\{ \tau \leq t \}$ is  measurable 
with respect to the  sigma field generated by $\left(\tilde{X}_{s} \right)_{0 \leq s \leq t} $.
Then, the problem considered in the following is to minimize
$\mathbb{E} \left[ R(\tau,X) \right] = \mathbb{E} \left[ I(\tau,X) + (1 -\tau )X \right]$, where $\tau < 1$ $a.s.$
Consider two arbitrary positive, strictly increasing and continuous functions $f_1$, $f_2$ on $[0, 1)$ with
$\int_{0}^{1} f_2(t) = \infty$. 
Let the random variable $Y$ be the height of the lowest $\mathscr{P}$-point above $f_1$ on $\left[ 0, \alpha \right]$, where $0 \leq \alpha<1$. That is: 
\begin{eqnarray}
\label{eq:ref1}
Y &:=& \min \left\{ X \ : \ \ (T,X) \in \mathscr{P}, \ \ T \in [0,\alpha],  \ \ X > f_1(T)  \right\}  \ .
\end{eqnarray}
\noindent In the following, $y$ stands for a realization of $Y$. Next we consider a stopping point $(\tau,X)$ with
\begin{eqnarray}
\label{eq:ref3}
\tau := \inf \left\{T \ : \ \ (T,X) \in \mathscr{P}, \ \ X { \leq } \left(f_1(T)  \wedge Y \right) \mathbf{1}_{(T \leq \alpha)} + f_2 (T) \mathbf{1}_{(T > \alpha)}   \right\}  \ ,
\end{eqnarray}
%
%
%
\noindent where the threshold function is random due to $Y$ . Note that $\tau < 1$ $a.s.$ by the condition $\int_0^1 f_2 = \infty$. 
By the stopping rule defined in (\ref{eq:ref3}), the relative rank of $X$ is $1$ in the event $\{ \tau \leq \alpha \}$.
Since $f_1$ is increasing, $y \leq f_1(t)$ implies that $(t',y) \in \mathcal{P}$, for some $0 \leq t' < t$. Hence,  $\tau$ is a stopping time for the filtration 
generated by $(\tilde X_s)_{0 \leq s \leq t}$. 
%

\begin{figure}[ht]
\begin{center}
\vskip-0.5em
\subfigure[Stopping Rule $\tau$. $f_1=\frac{b}{1-t}$, $f_2=\frac{c}{1-t}$.]{\includegraphics[height=57mm]{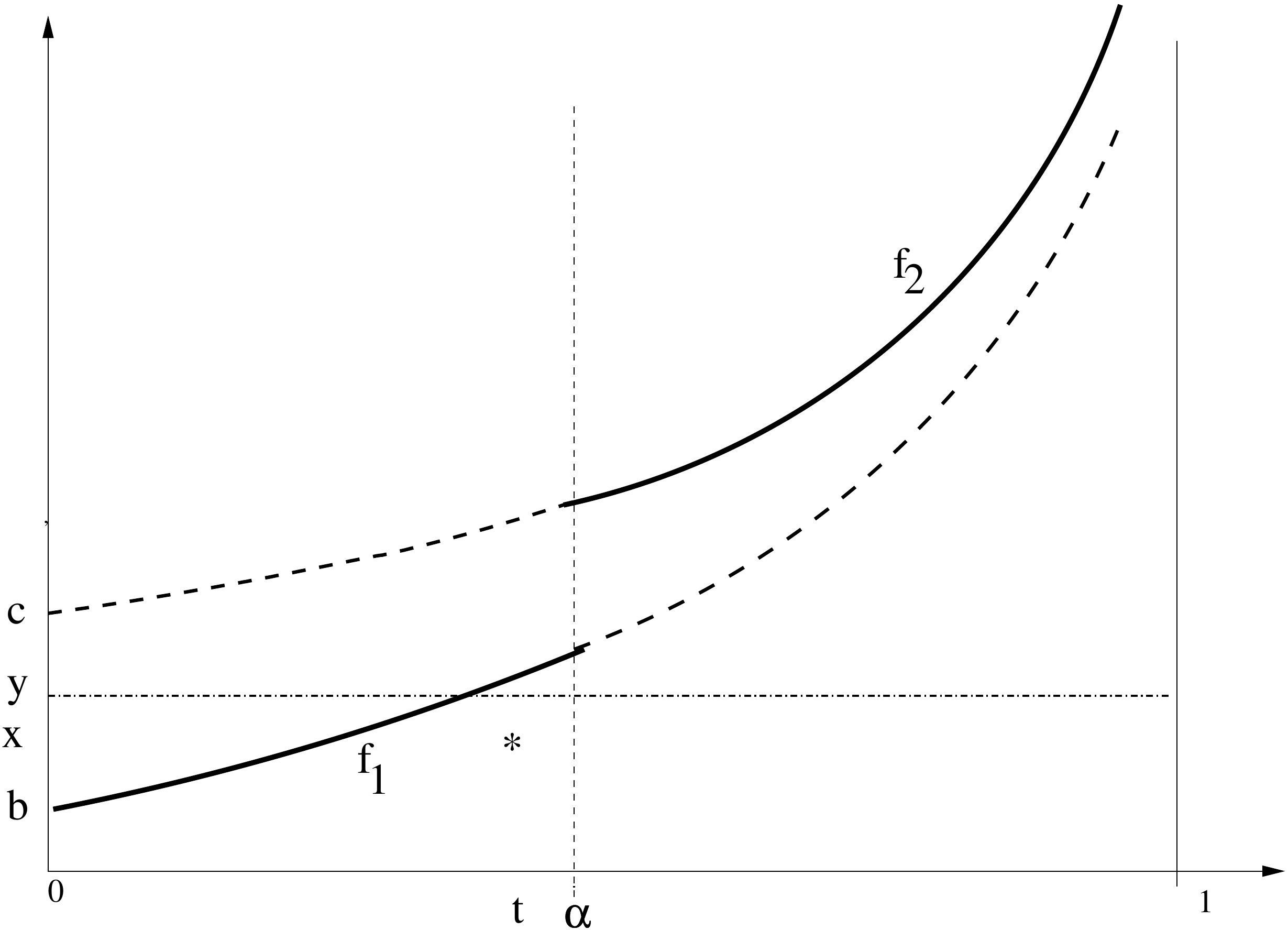} \label{subfig10} }
\subfigure[Area between $y$ and and the graph of $f_1 (t)$.]{\includegraphics[height=57mm]{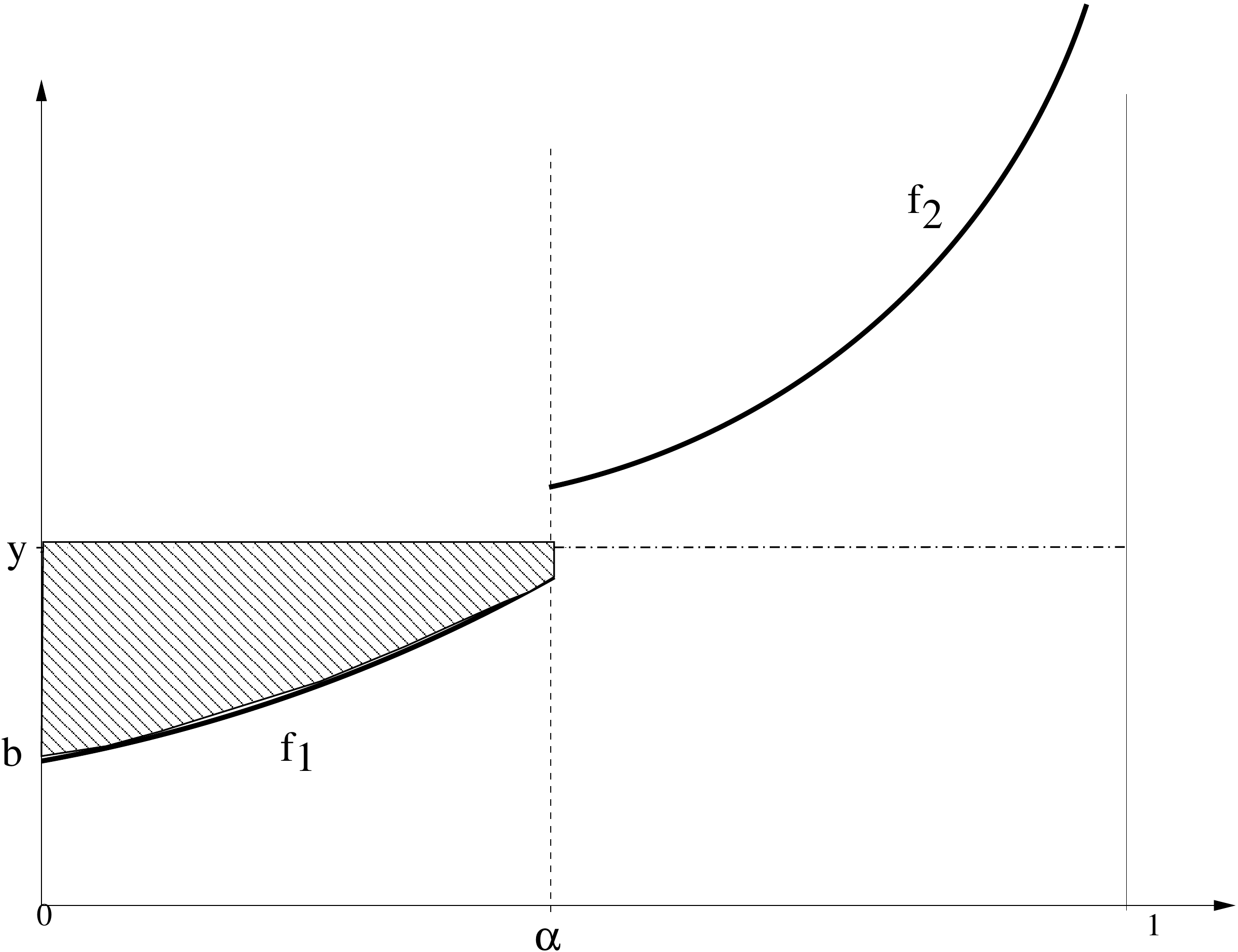} 
 \label{figsub12a} } \\[1pt]
\subfigure[Loss components, $y > f_2(t)$ is a realizations of $Y$.]{\includegraphics[height=57mm]{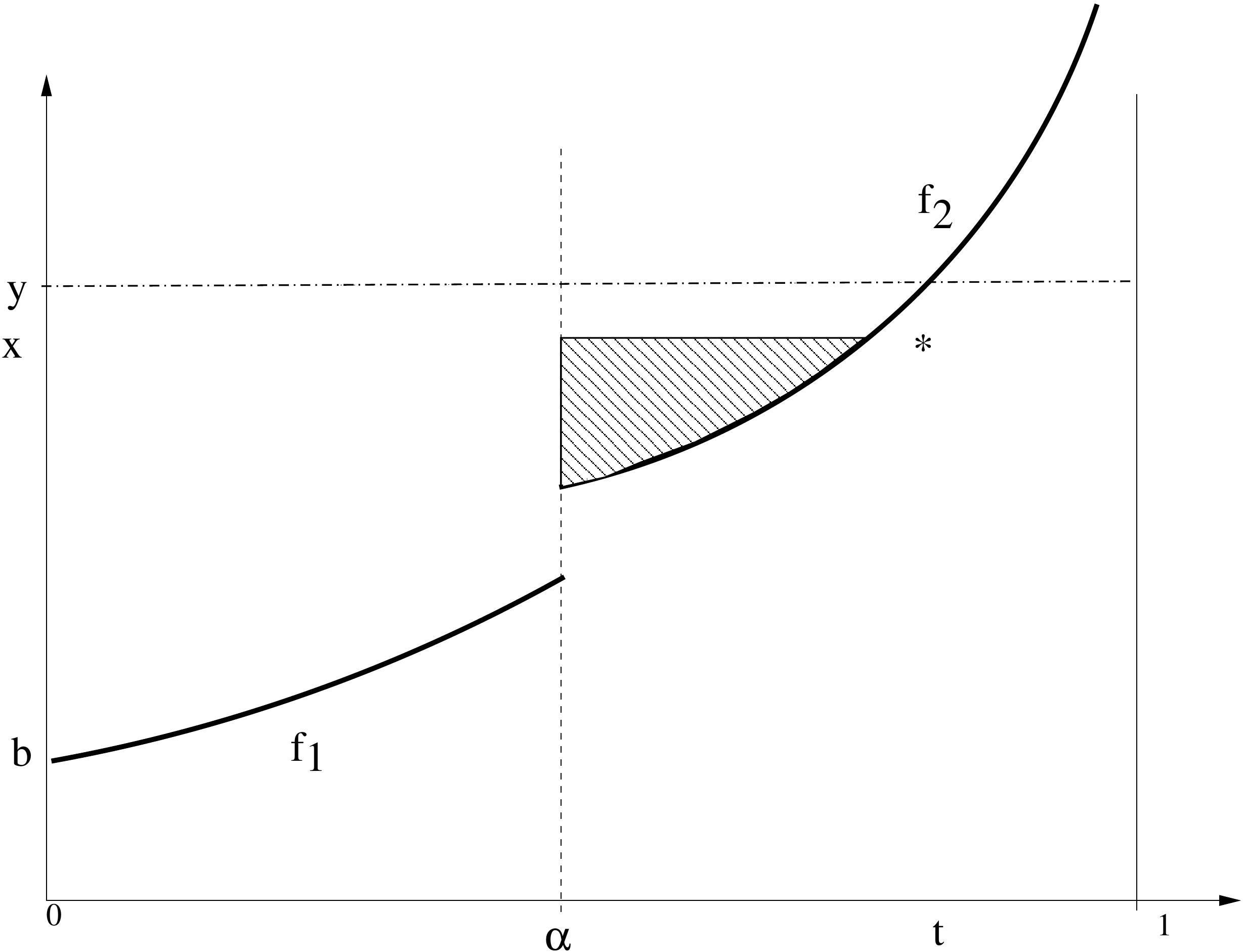} \label{figsub17a} }
\subfigure[Loss components, $y < f_2(t)$ is a realiations of $Y$.]{\includegraphics[height=57mm]{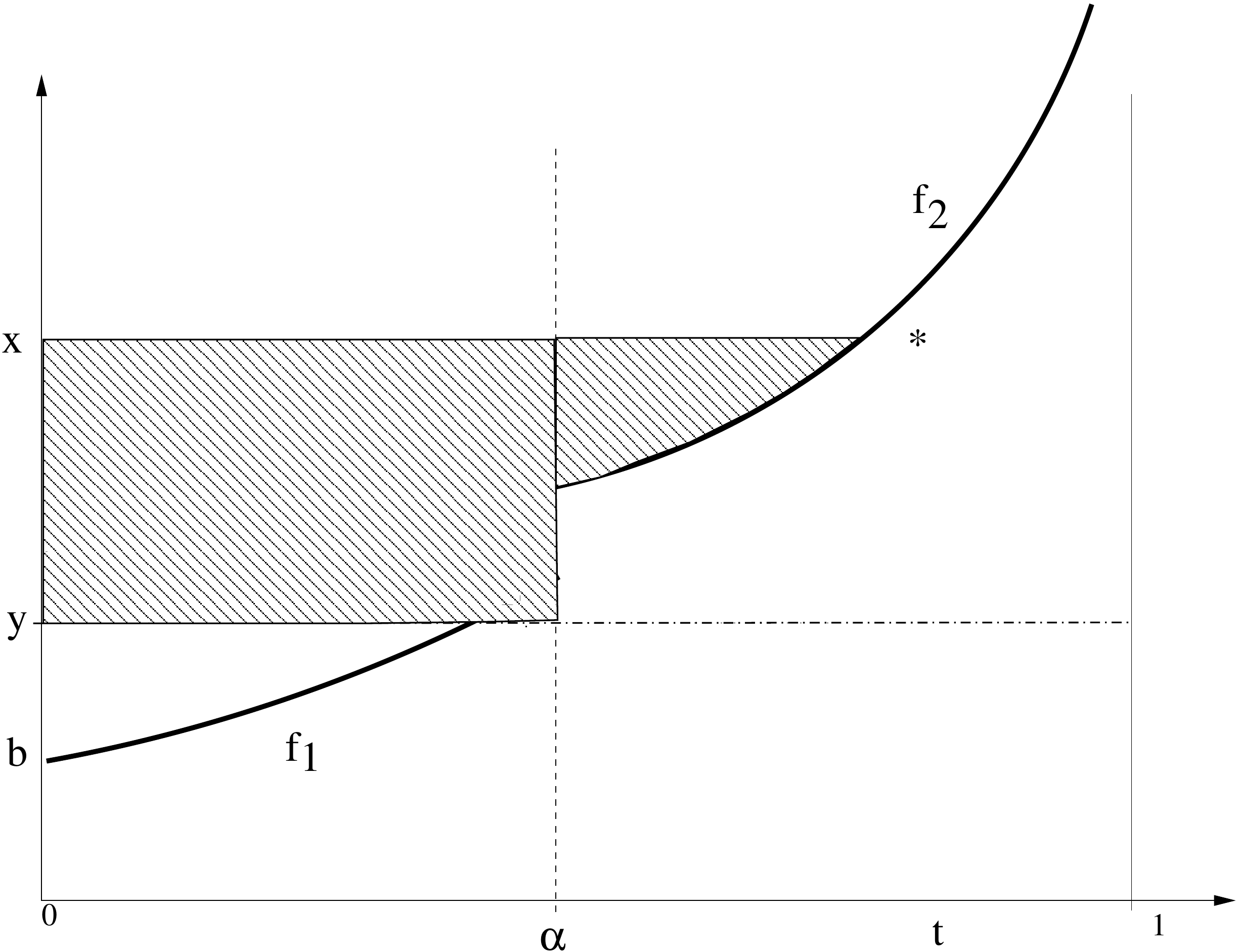} \label{figsub17b} }
\vskip-1em
\vskip2em
\caption{The stopping rule $\tau$ and loss components.}
\end{center}
\end{figure}

Figure~\ref{subfig10} provides a graphical description of the stopping rule $\tau$. For $t \leq \alpha$, a decision maker stops if a particle $(t,x)$ with $x \leq y$ and
$x \leq f_1(t)$ is observed, while for $ t > \alpha$, the threshold $f_2$ is applied. Figure~\ref{subfig10} shows one realization of $Y$, where $y<f_1(\alpha)$. 
The decision whether to stop with the first $(t,x)$ such that $x \leq f_1(t)$ and $t \leq \alpha$ depends on $y$ only if
$t \in (f_1^{-1} (y)  ,\alpha]$.
Figure~\ref{figsub12a} provides an example where $y>f_1(t)$
for all $t \leq \alpha$. In this case, the decision maker stops with the first $(t,x)$, such that $0 \leq t \leq \alpha$ and $x \leq f_1(t)$, if any such $(t,x)$ exists. 

To obtain the risk $\mathbb{E} \left[ R( \tau,X) \right]$, we consider the conditional risk $\mathbb{E} \left[ R( \tau,X) |Y = y \right]$ and integrate out $y$.
To do this, we define
\begin{eqnarray}
\label{eq:ref5}
F_1( t,y) &:=& \int_{0}^{t} \left(f_1(s)  \wedge y \right) ds \ , \ \ \ 
F_2 ( t) := \int_{  \alpha }^{t} f_2(s)   ds \ , \text{ and } \\ 
\label{eq:ref7}
S_1( x) &:=& \int_{0}^{x} \left(f_1^{-1}(z)  \wedge \alpha \right) dz = \alpha x - F_{1} (\alpha,x) \ , \ 
\end{eqnarray}
\noindent where $f_1^{-1}(z)$ = 0 for $z < f_1(0)$. Furthermore, 
\begin{eqnarray}
\label{eq:ref8}
S_2(x)  &:=& \int_{0}^{x} \left(f_2^{-1}(z) - \alpha \right) dz  \ , \ 
\end{eqnarray}
\noindent where $f_2^{-1}(z) =  \alpha$ for $z < f_2( \alpha)$. 
By the properties of the planar Poisson process the probability that the area bounded by the horizontal line with  
height $y$ and the graph of $f_1$ for the time interval $[0,\alpha]$ is empty is $e^{- S_1(y)}$; see the shaded area in Figure~\ref{figsub12a}.
Hence, the density of $Y$ is 
\begin{eqnarray}
\label{eq:ref9}
\mathbb{P} \left( Y \in dy \right)  &=&  \left[ \frac{d}{dy} \left( 1 - e^{- S_1(y)} \right) \right] dy = e^{- S_1(y)} \left(f_1^{-1}(y) \wedge \alpha \right)_{} dy  \ . \ 
\end{eqnarray}
By the definition of $Y$, $\mathbb{P} \left( Y \in dy \right)>0$, for $y > f_1(0)$ and zero else.  
Given $\{ Y =y \}$, the conditional joint density of $(\tau,X)$ is
%
\begin{eqnarray}
\label{eq:ref9}
\mathbb{P} \left( (\tau,X)  \in (dt,dx)| Y=y \right)  &=& 
\begin{cases}
e^{- F_1(t,y)} dt \ dx &  0 \leq t \leq \alpha,  \ \ \ 0 \leq x \leq f_1^{}(t) \wedge y \ , \\ 
e^{- F_1(\alpha,y) - F_2(t)} dt \ dx &   \alpha < t \leq 1,  \ \  \ 0 \leq x \leq f_2^{}(t)  \ . \
\end{cases}	
\end{eqnarray} 
\noindent By means of the density (\ref{eq:ref9}), we obtain the conditional probability $\mathbb{P} \left( \tau > t | Y=y \right) = e^{- F_1(t,y)}$, for $t \leq \alpha$, 
and $\mathbb{P} \left( \tau > t | Y=y \right) = e^{- F_1(\alpha,y) - F_2(t)}$, for $ \alpha< t \leq 1$. 
Splitting the time at $\alpha$, the two components of the risk given $Y = y$ are computed as
%
	%
	\begin{eqnarray}
	\label{eq:ref10}
	\mathbb{E} \left[ R (\tau,X) \mathbf{1}_{(\tau \leq \alpha)}  | Y=y \right]  &=& 
	\int_{0}^{\alpha}
	e^{- F_1(t,y)} 
	\left(
	\int_{0}^{f_1(t) \wedge y} 
	\left(1 + x(1-t) \right) \ dx \right) \ dt 
	\ \text{ and } \nonumber \\
	\mathbb{E} \left[ R (\tau,X) \mathbf{1}_{(\tau > \alpha)}  | Y=y \right]  &=& 
	\int_{\alpha}^{1}
	e^{- F_1(\alpha,y) - F_2(t)}  
	\left[
	\int_{0}^{f_2(t) \wedge y} 
	\left(1 + x(1-t)  \right) \ dx
	+
	\int_{f_2(\alpha)}^{f_2(t)} 
	S_2 (x)   \ dx
    \right. \nonumber \\ && + \left.
	\int_{f_2(t) \wedge y}^{f_2(t)} 
	\left( 2 + x (1-t) +  
	(x -y ) \alpha  \right) \ dx 
	\right]  \ dt \ .
	\end{eqnarray}%
%
For $ t \leq \alpha$, the integrand of the inner integral describes the conditional expected loss contribution when we stop at $(t,x)$. In this case, the relative rank of $(t,x)$ is $1$ and the conditional expected loss of $(t,x)$ is $1$ plus the expected number of atoms in the south-east of $(t,x)$. That is, $1 + x(1-t)$.  An example for such a $(t,x)$ is provided by ``$*$'' in  Figure~\ref{subfig10}.
 For $t \geq \alpha$, we have to account for atoms in the area above the curve
$f_2(s)$, $ \alpha < s \leq t$ and below the level obtained by $x$. In the Figures~\ref{figsub17a} 
and \ref{figsub17b} this corresponds to the almost triangular shaded area. Formally, the size of this area is 
$S_2 (x)$, where $S_2(x)=0$, for $x \leq f_2(\alpha)$ and $S_2(x)>0$, for $x > f_2(\alpha)$. In particular, 
Figure~\ref{figsub17a} shows $S_2 (x)$, for some $x > f_2(\alpha)$.  
Figures~\ref{figsub17a} and \ref{figsub17b} describe two generic cases where a decision maker stops at $(t,x)$ with $t > \alpha$.  
Figure~\ref{figsub17a} describes the case $y>x$. The loss coming from the past is described by the shaded almost triangular area. The conditional expected loss  
is $1+ x(1-t) + S_2(x)$. $1$ accounts for the particle itself and $x(1-t)$ for the expected loss in the future, which is the size of the area in the south-east of the particle $(t,x)$. 
Next, Figure~\ref{figsub17b} describes the case where $y \leq x$. In this case, the loss coming from the past is described by the shaded almost triangular area
and the shaded rectangular area. The conditional expected loss  
is $2+ x(1-t) + S_2(x) + (x-y) \alpha$. $2$ accounts for the particle itself and the particle $(t',y)$, such that $y$ is the realization of $Y$.  In addition, $x(1-t)$ accounts for the expected loss in the future.
Last but not least, $F_1( t,y) $ and $F_2(t)$ measure the size of areas below $f_1(s) \wedge y$, where $0 \leq  s \leq t$,
and $f_2(s)$, where $\alpha<s \leq t$. $S_1(y)$ measures to the difference of the size of the area of the 
rectangle described by the points $(0,0)$, $(\alpha,0)$, $(\alpha,f_1(\alpha) \wedge y )$ and $(0,f_1(\alpha) \wedge y )$ and the size of the area
below the graph $\{ \left( t,f_1(t) \wedge y \right) : \ 0 \leq t \leq \alpha \}$. The shaded area in Figure~\ref{figsub12a} provides an example for this area. 
  
Next, we work with $f_1(t) = \frac{b}{1 - t}$ and $f_2(t) = \frac{c}{1 - t}$. For these functions, 
we were able to obtain the conditional expected losses $
\mathbb{E} \left[ R (\tau,X) \mathbf{1}_{(\tau \leq \alpha)}  | Y=y \right]$ and  
$\mathbb{E} \left[R (\tau,X) \mathbf{1}_{(\tau > \alpha)}  | Y=y \right]$, as well as the expected loss 
$\mathbb{E} \left[ R (\tau,X) \right] = \int_b^{\infty}  
\left( \mathbb{E} \left[ R (\tau,X) \mathbf{1}_{(\tau \leq \alpha)}  | Y=y \right]+
\mathbb{E} \left[ R (\tau,X) \mathbf{1}_{(\tau \leq \alpha)}  | Y=y \right] \right)
e^{- S_1(y)} \left(f_1^{-1}(y) \wedge \alpha \right)_{} dy $ in closed form by using the \texttt{Mathematica~8.0} package. 
%
%
%
%
 %
%
%
%
%

To minimize $\mathbb{E} \left[ R (\tau,X) \right]$ the parameters $\alpha$, $b$ and $c$ have to be chosen optimally.%
{\footnote{
To do this, the {\tt MATHEMATICA} expressions were converted to {\tt MATLAB R2012a} code.   }} %
By numerical tools we observed that the value of the loss function is approximately minimized with $(\alpha_{*},b_{*},c_{*})^{} = (
0.34328,
1.82571,
2.00000
)^{}$. With these parameters we calculated an expected loss of $2.32614$.

\begin{remark}
To check the above results we performed various simulation studies. We observe that our rule with
$(\alpha_{*},b_{*},c_{*})^{} = (
0.34328,
1.82571,
2.00000
)^{}$ dominates the rule of \cite{Tamaki2004}, where 
$(\alpha_{T},b_{T},c_{T})^{} = (
0.42,
1.95,
1.95
)^{}$, as well as the threshold rule presented in \cite{Gnedin2007}, where $(\alpha_{G},b,c_{G})^{} = (
0,
b,
1.9469
)^{}$ and $b>0$ is arbitrary, in the mean. For this threshold rule the expected rank is $1+ \frac{c_G}{2} + \frac{1}{c_G^2-1}=
2.3318$. 

At the end of their paper, \cite{AssafChan1996} mentioned in Remark 6.2  that they tried a rule that would be in retrospect the $n$-period analogue of our rule, but with $b=c$
(as already mentioned, this rule has been investigated by \cite{Tamaki2004}, where a value of $2.33044$, with $\alpha_T=0.42$ and $b_T=c_T=1.95$, has been obtained in closed form). That is, they tried thresholds, where in the beginning fraction of time $\gamma$ the additional condition of stopping only with relative rank $1$ is imposed, and then a threshold rule is used for the remaining time. However, they left the thresholds at $\min \{ 2/(n-k+1), 1 \}$, for $k=1,\dots,n$ (these thresholds had been known to be good approximations for thresholds of the optimal threshold rule for large enough $n$).

\cite{AssafChan1996} reported that ``the improvement is however very small, and even with $10,000$ simulations the
standard error is too large to determine whether the improvement is real.''  Also in our simulation runs (with 50,000 steps), we observed that the standard errors are large compared to the improvements  
obtained with our rule and the rule of \cite{Tamaki2004}, which also implies that a reliable comparison of rules can hardly be performed by means of a simulation analysis.

In addition, we inserted the parameters 
$(\alpha_{T},b_{T},c_{T})^{}$ and  $(\alpha_{G},b,c_{G})^{}$ into our closed form expression obtained for $\mathbb{E} \left[ R (\tau,X) \right]$. 
For the threshold rule, the difference  between the number we get by inserting the corresponding parameters into the expression obtained for
$\mathbb{E} \left[ R (\tau,X) \right]$ and the number we get by means of  
$1+ \frac{c_G}{2} + \frac{1}{c_G^2-1}$ is smaller than $4 \cdot 10^{-8}$. For 
$(\alpha_{T},b_{T},c_{T})^{}$ we obtained $2.33045$ instead of $2.33044$ derived by \cite{Tamaki2004}. That is, up to a numerical error, 
 the expected losses obtained in \cite{Tamaki2004} and \cite{Gnedin2007} can also be replicated by means of the closed form expression in this article.
\end{remark}

\clearpage
\bibliographystyle{apalike}
\bibliography{master} 

\clearpage

\appendix
\section{Mathematica Code}
The following pages present the Mathematica code used to obtain the expected loss $\mathbb{E} \left[ R (\tau,X) \right]$:



\section*{Supplementary material}

```mathematica
(*Calculations to "A new strategy for Robbins' problem of optimal stopping"
 by Martin Meier and Leopold Sögner, April 2016*)

(* start part 1, some definitions *)
(* Obtain f1, f2 F1(t,y), F2(t), S1(x) and S2(x) *)
f1[s] = b / (1 - s) ; (* fb(t) *)

f2[s] = c / (1 - s) ;

F11 = Integrate[f1[s], {s, 0, t}, Assumptions →
    {b > 0, s ∈ Reals, t ∈ Reals, y ∈ Reals, 0 < t ≤ 1, b > 1, y > 0, y > b / (1 - t)}]
 (* $ F_1(t,y) $, defined in equation (3), for $f_1 \leq y$ *)
```

$-b \, \text{Log}[1-t]$

```mathematica
F11alphay = Integrate[f1[s], {s, 0, alpha}, Assumptions →
    {b > 0, s ∈ Reals, t ∈ Reals, y ∈ Reals, 0 < t ≤ 1, b > 1, y > 0, y > b / (1 - t), alpha < 1}]
 (* $ F_1(\alpha,y) $, defined in equation (3), for $f_1 \leq y$ *)
```

$-b \, \text{Log}[1-\text{alpha}]$

```mathematica
F12 = b - y + t y - b Log[b / y] (* $ F_1(t,y) $,
defined in equation (3), for $f_1(\alpha) > y$ *)
F12alphay = b - y + alpha y -
    b Log[b / y] (* $ F_1(\alpha,y) $, defined in equation (3), for $f_1 > y$ *)
```

$b - y + t\, y - b \, \text{Log}\!\left[\dfrac{b}{y}\right]$

$b - y + \text{alpha}\, y - b \, \text{Log}\!\left[\dfrac{b}{y}\right]$

```mathematica
F2 = Integrate[c / (1 - s), {s, alpha, t}, Assumptions → {c > 0, alpha > 0, 1 > t > alpha > 0}]
 (* $ F_2(t) $, defined in equation (3)*)
```

$c \, \text{Log}\!\left[\dfrac{-1+\text{alpha}}{-1+t}\right]$

```mathematica
(*Define Case I: $y \leq b/(1-\alpha) =
 f_1(\alpha) $ and Case II: $y > b/(1-\alpha) = f_1(\alpha) $ *)

S1caseI = alpha y + b Log[b / y]
(* $S_1 (x) $ defined in equation (4) evaluated at $x=y$,
for the case $y \leq b/(1-alpha)$ *)
S1caseII = Simplify[alpha y - (b - y + alpha * y - b * Log[b / y])]
(* $S_1 (x) $ defined in equation (4) evaluated at $x=y$,
for the case $y > b/(1-alpha)$ *)
```

$\text{alpha}\, y + b \, \text{Log}\!\left[\dfrac{b}{y}\right]$

$-b + y + b \, \text{Log}\!\left[\dfrac{b}{y}\right]$

```
f2inv[z] = 1 - c / z;
S2[x] =
 Piecewise[{{Integrate[(f2inv[z] - alpha), {z, c / (1 - alpha), x}, Assumptions → {c > 0,
       x > c > 0, alpha > 0, c + alpha x < x}], x > c / (1 - alpha)}, {0, x < c / (1 - alpha)}}]
(* $f_2^{-1}$ and $S_2 (x) $ defined in equation (5) *)
```

$$\begin{cases} -c + x - \text{alpha}\, x - c\, \text{Log}\!\left[\dfrac{x-\text{alpha}\, x}{c}\right] & x > \dfrac{c}{1-\text{alpha}} \\ 0 & \text{True} \end{cases}$$

```
(* end first step*)

(* Second step Obtain Conditional expected loss (conditional on $\{ Y=y \}$ ) *)


(* inner integral for the $t \leq alpha$ case;
this  is the first integral in equation (8);
terms referring to this first integral contain I1 *)
I1integrad = 1 + x (1 - t);

I1integralc1 =
 Simplify[Integrate[I1integrad, {x, 0, b / (1 - t)}, Assumptions → {x > 0, 0 < t < 1}]]
  (* Inner Integral for the $f_1(t) \leq y$ case *)
```

$$-\dfrac{b\,(2+b)}{2\,(-1+t)}$$

```
I1integralc2 = Simplify[Integrate[I1integrad, {x, 0, y}, Assumptions → {x > 0, 0 < t < 1}]]
(* Inner Integral for the $f1(t) > y$ case *)
```

$$\dfrac{1}{2}\, y\,(2 + y - t\, y)$$

```
Innerintegral1 = Piecewise[ {{ I1integralc1, t ≤ 1 - b / y}, { I1integralc2, t > 1 - b / y} } ]
```

$$\begin{cases} -\dfrac{b\,(2+b)}{2\,(-1+t)} & t \leq 1 - \dfrac{b}{y} \\ \dfrac{1}{2}\, y\,(2 + y - t\, y) & t > 1 - \dfrac{b}{y} \\ 0 & \text{True} \end{cases}$$

```
(* outer integral , $t \leq  \alpha$ *)
(* outer integral , $t \leq alpha$, $y \leq b/(1-alpha)$ and $f1 \leq < y$;
then $f_1 \leq y$ for $t \leq 1-b/y$ *)

Outer1integralcaseI1a =  Simplify[Integrate[Exp[-F11] * I1integralc1,
    {t, 0, 1 - b / y}, Assumptions → {x > 0, 0 < t < alpha < 1, 0 < b < y}]]
```

$$-\dfrac{1}{2}\,(2+b)\left(-1 + \left(\dfrac{b}{y}\right)^{b}\right)$$

```
(* outer integral , $t \leq \alpha$ *)(* outer integral ,
$t \leq alpha$, $y \leq b/(1-alpha)$ for $t >1-b/y$ we have $y < f_1$ *)

Outer1integralcaseI1b =
 Simplify[Integrate[Exp[-F12] * I1integralc2, {t, 1 - b / y, alpha},
    Assumptions → {x > 0, 0 < t < alpha < 1, 0 < b < y < b / (1 - alpha)}]]
```

$$\dfrac{1}{2}\left(\dfrac{b}{y}\right)^{b}\left(1 + b + e^{-b+y-\text{alpha}\, y}\,(-1 + (-1 + \text{alpha})\, y)\right)$$

```
(* outer integral , $t \leq \alpha$ *)(* outer integral , $t \leq alpha$,
$y > b/(1-alpha)$, i.e. $f_1 \leq y$ for all $t \in [0,\alpha]$ *)
```

```
Outer1integralcaseII24 = Simplify[Integrate[Exp[-F11] * I1integralcase1,
   {t, 0, alpha}, Assumptions → {x > 0, 0 < t < alpha < 1, b / (1 - alpha) < y}]]
(* f1(t) < y and y > b/(1-alpha) *)
```

$$\frac{\left(1 - (1 - \text{alpha})^{1+b}\right) \text{I1integralcase1}}{1 + b}$$

```
(* end: outer integral , $t \leq \alpha$ case *)

(* begin: integrals , $t > alpha$ *) (* first term in equation (8)→ 21 *)
I21integralcase1 =
 Simplify[Integrate[I1integrad, {x, 0, c / (1 - t)}, Assumptions → {x > 0, 0 < t < 1}]]
(* $f_2(t) \leq y$ *)
I21integralcase2 =
 Simplify[Integrate[I1integrad, {x, 0, y}, Assumptions → {x > 0, 0 < t < 1}]]
(* $f_2(t)  > y$ *)
```

$$-\frac{c(2+c)}{2(-1+t)}$$

$$\frac{1}{2} y (2 + y - t y)$$

```
(* second term in equation (8) → 22; i.e. we obtain integrals of $ S_2(x) $ *)

I22integral = Simplify[Integrate[S2[x],
   {x, 0, c / (1 - t)}, Assumptions → {x > 0, c > 0, 0 < t < 1, 0 < alpha < t}]]
```

$$-\left(c^2\left(-2\,\text{alpha} + \text{alpha}^2 - (-2+t)\,t + 2(-1+\text{alpha}+t-\text{alpha}\,t)\,\text{Log}\left[\frac{-1+\text{alpha}}{-1+t}\right]\right)\right) / \left(2(-1+\text{alpha})(-1+t)^2\right)$$

```
I22integralcheck = Simplify[Integrate[S2[x], {x, c / (1 - alpha), c / (1 - t)},
   Assumptions → {x > 0, c > 0, 0 < t < 1, alpha > 0, alpha < t}]]
```

$$-\left(c^2\left(-2\,\text{alpha} + \text{alpha}^2 - (-2+t)\,t + 2(-1+\text{alpha}+t-\text{alpha}\,t)\,\text{Log}\left[\frac{-1+\text{alpha}}{-1+t}\right]\right)\right) / \left(2(-1+\text{alpha})(-1+t)^2\right)$$

```
(* third term in equation (8)→ 23;
we have to tell Mathematica that this term is >0 if f2(t)>y and zero 0 else *)

I23integrand = 2 + x (1 - t) + (x - y) alpha
```

$$2 + (1-t)\,x + \text{alpha}\,(x-y)$$

```
I23integralf2smallery = 0 ;(* f2(t) \leq y *)
I23integralf2largery =
 Simplify[Integrate[I23integrand, {x, y, c / (1 - t)}, Assumptions → {x > 0, 0 < t < 1}]]
(* f2(t) > y *)
```

$$\frac{(c + (-1+t)\,y)\,(c\,(1+\text{alpha}-t) + (-1+t)\,(-4 + (-1+\text{alpha}+t)\,y))}{2(-1+t)^2}$$

```
(* end inner second integrals*)

(* begin outer second integrals, i.e when $t > \alpha$ *)
```

```
(* In the following we consider 4 cases:
   Case 1: $f_1(0) =
 b \leq y \leq c/(1-\alpha) = f_2(\alpha)$ and $y \leq b/(1-\alpha) = f_1(\alpha)$
      Case 2: $f_1(0) =
    b \leq y \leq c/(1-\alpha) = f_2(\alpha)$ and $y > b/(1-\alpha) = f_1(\alpha)$
           Case 3: $f_1(0) =  y > c/(1-\alpha) =
         f_2(\alpha)$ and $y \leq b/(1-\alpha) = f_1(\alpha)$
              Case 4: $f_1(0) = y > c/(1-\alpha) =
            f_2(\alpha)$ and $y > b/(1-\alpha) = f_1(\alpha)$
             In addition: we consider the case where $c$ is very small,
i.e. $c/(1-\alpha) < b$
*)


(*In the following integrals we use
 F12alphay the cases 1 and 3 ($y \leq b/(1-\alpha)$),
while F11alphay has to be used in the cases 2 and 4 ($y > b/(1-\alpha)$);
depending on whether $y \leq c/(1-\alpha)$ (cases 1 and 2) or $y > c/(1-\alpha)
  $ (cases 3 and 4) we use the corresponding integrands obtained above    *)

Outer21integralcase1 =
 Simplify[Exp[-F12alphay] * Integrate[Exp[-F2] * I21integralcase2, {t, alpha, 1},
    Assumptions → {0 < alpha < t < 1, 0 < c, y > b > 0, y < c / (1 - alpha)}]]
```

$$\frac{1}{2} \, (-1 + \text{alpha}) \, e^{-b + y - \text{alpha}\, y} \left(\frac{b}{y}\right)^b y \left(-\frac{2}{1+c} + \frac{(-1+\text{alpha})\, y}{2+c}\right)$$

```
Outer21integralcase1csmall =
 Simplify[Exp[-F12alphay] * Integrate[Exp[-F2] * I21integralcase1, {t, alpha, 1 - c / y},
    Assumptions → {0 < alpha < t < 1, 0 < c, y > b > 0, y < c / (1 - alpha)}]] +
  Simplify[Exp[-F12alphay] * Integrate[Exp[-F2] * I21integralcase2, {t, 1 - c / y, 1},
    Assumptions → {0 < alpha < t < 1, 0 < c, y > b > c > 0, y < c / (1 - alpha)}]]
```

$$\frac{c\,(4 + 3c + c^2)\, e^{-b+y-\text{alpha}\,y} \left(\frac{b}{y}\right)^b \left(\frac{c}{y - \text{alpha}\,y}\right)^c}{2\,(1+c)\,(2+c)} - \frac{1}{2}(2+c)\, e^{-b+y-\text{alpha}\,y} \left(\frac{b}{y}\right)^b \left(-1 + \left(\frac{c}{y - \text{alpha}\,y}\right)^c\right)$$

```
Outer21integralcase2 =
 Simplify[Exp[-F11alphay] * Integrate[Exp[-F2] * I21integralcase1, {t, alpha, 1 - c / y},
    Assumptions → {0 < alpha < t < 1, 0 < c, y > b > 0, y < c / (1 - alpha)}]] +
  Simplify[Exp[-F11alphay] * Integrate[Exp[-F2] * I21integralcase2, {t, 1 - c / y, 1},
    Assumptions → {0 < alpha < t < 1, 0 < c, y > b > c > 0, y < c / (1 - alpha)}]]
(* note that for $t \in [\alpha,1-c/y]$ we have $y \geq f_2(t)$,
    while for $t \in (1-c/y,1]$ we observe that $y < f_2(t)$ *)
```

$$\frac{(1-\text{alpha})^b \, c\, (4 + c\,(3+c)) \left(\frac{c}{y - \text{alpha}\,y}\right)^c}{2\,(1+c)\,(2+c)} - \frac{1}{2}(1-\text{alpha})^b (2+c) \left(-1 + \left(\frac{c}{y - \text{alpha}\,y}\right)^c\right)$$

```
Outer21integralcase3 = Simplify[
   Simplify[Exp[-F12alphay] * Integrate[Exp[-F2] * I21integralcase1, {t, alpha, 1 - c / y},
       Assumptions → {x > 0, 0 < alpha < t < 1, 0 < c, y > c / (1 - alpha), y > b > 0}]] +
    Simplify[Exp[-F12alphay] * Integrate[Exp[-F2] * I21integralcase2, {t, 1 - c / y, 1},
       Assumptions → {x > 0, 0 < alpha < t < 1, 0 < c, y > c / (1 - alpha), y > b > 0}]]]
```

$$\frac{1}{2\,(1+c)\,(2+c)}$$

$$e^{-b+y-\text{alpha}\,y}\left(\frac{b}{y}\right)^b\left(4+c^3-4\left(\frac{c}{y-\text{alpha}\,y}\right)^c+c^2\left(5-2\left(\frac{c}{y-\text{alpha}\,y}\right)^c\right)-4\,c\left(-2+\left(\frac{c}{y-\text{alpha}\,y}\right)^c\right)\right)$$

```
Outer21integralcase4 = Simplify[
   Simplify[Exp[-F11alphay] * Integrate[Exp[-F2] * I21integralcase1, {t, alpha, 1 - c / y},
       Assumptions → {x > 0, 0 < alpha < t < 1, 0 < c, y > b > 0, y > c / (1 - alpha)}]] +
    Simplify[Exp[-F11alphay] * Integrate[Exp[-F2] * I21integralcase2, {t, 1 - c / y, 1},
       Assumptions → {x > 0, 0 < alpha < t < 1, 0 < c, y > b > 0, y > c / (1 - alpha)}]]]
```

$$\frac{(1-\text{alpha})^b\left(\frac{c\,(4+c\,(3+c))\left(\frac{c}{y-\text{alpha}\,y}\right)^c}{1+c}-(2+c)^2\left(-1+\left(\frac{c}{y-\text{alpha}\,y}\right)^c\right)\right)}{2\,(2+c)}$$

(* Since the inner integrals for the 22 terms (i.e. $\int S_2(x)\ dx$)
 do not depend on $y$ we can consider the cases 1 and 3 as and
 well as 2 and 4 joints to obtain the integral for the 22 term*)

```
Outer22integralcase24 =
  Simplify[Exp[-F11alphay] * Integrate[Exp[-F2] * I22integral, {t, alpha, 1},
      Assumptions → {0 < alpha < t < 1, b / (1 - alpha) < y, c > 1, b > 0}]]
```

$$\frac{(1-\text{alpha})^b}{-1+c^2}$$

```
Outer22integralcase13 =
  Simplify[Exp[-F12alphay] * Integrate[Exp[-F2] * I22integral, {t, alpha, 1},
      Assumptions → {0 < alpha < t < 1, b / (1 - alpha) > y > b, c > 1, b > 0}]]
```
**for**
 **t**

$$\frac{e^{-b+y-\text{alpha}\,y}\left(\frac{b}{y}\right)^b}{-1+c^2}$$

for t

```
Outer23integralcase1 =
  Simplify[Exp[-F12alphay] * Integrate[Exp[-F2] * I23integralf2largery, {t, alpha, 1},
      Assumptions → {0 < alpha < t < 1, b / (1 - alpha) > y > b, c > 1, b > 0}]]
```

$$\frac{1}{2}\,e^{-b+y-\text{alpha}\,y}\left(\frac{b}{y}\right)^b\left(\frac{-(-1+c)\,(4+c)+\text{alpha}\,(-4+3\,c)}{(-1+\text{alpha})\,(-1+c)}+\right.$$
$$\left.\frac{2\,(-2+\text{alpha}-\text{alpha}\,c)\,y}{1+c}-\frac{(-1+\text{alpha})\,(-1-c+\text{alpha}\,(3+2\,c))\,y^2}{(1+c)\,(2+c)}\right)$$

```
Outer23integralcase1csmall =
 Simplify[Exp[-F12alphay] * Integrate[Exp[-F2] * I23integralf2largery,
    {t, 1 - c / y, 1}, Assumptions → {0 < alpha < t < 1, 1 < c < y, y < c / (1 - alpha)}]]
```

$$\frac{e^{-b+y-\text{alpha}\,y}\left(\frac{b}{y}\right)^{b}\left(\frac{c}{y-\text{alpha}\,y}\right)^{c}\left(-4+c+2c^{2}+c^{3}+\text{alpha}\,(2+c)\,y\right)}{(-1+c)\,(1+c)\,(2+c)}$$

```
Outer23integralcase2 =
 Simplify[Exp[-F11alphay] * Integrate[Exp[-F2] * I23integralf2largery,
    {t, alpha, 1}, Assumptions → {0 < alpha < t < 1, b / (1 - alpha) < y, c > 1, b > 0}]]
```

$$\frac{1}{2}(1-\text{alpha})^{b}\left(\frac{-(-1+c)\,(4+c)+\text{alpha}\,(-4+3c)}{(-1+\text{alpha})\,(-1+c)}+\right.$$
$$\left.\frac{2\,(-2+\text{alpha}-\text{alpha}\,c)\,y}{1+c}-\frac{(-1+\text{alpha})\,(-1-c+\text{alpha}\,(3+2c))\,y^{2}}{(1+c)\,(2+c)}\right)$$

```
Outer23integralcase3 = Simplify[Exp[-F12alphay] *
   Integrate[Exp[-F2] * I23integralf2largery, {t, 1 - c / y, 1}, Assumptions →
    {0 < alpha < t < 1, b / (1 - alpha) > y > b, c > 1, b > 1, y > c / (1 - alpha)}]]
```

$$\frac{e^{-b+y-\text{alpha}\,y}\left(\frac{b}{y}\right)^{b}\left(\frac{c}{y-\text{alpha}\,y}\right)^{c}\left(-4+c+2c^{2}+c^{3}+\text{alpha}\,(2+c)\,y\right)}{(-1+c)\,(1+c)\,(2+c)}$$

```
Outer23integralcase4 =
 Simplify[Exp[-F11alphay] * Integrate[Exp[-F2] * I23integralf2largery, {t, 1 - c / y, 1},
    Assumptions → {0 < alpha < t < 1, b / (1 - alpha) < y, c > 1, b > 1, y > c / (1 - alpha)}]]
```

$$\frac{(1-\text{alpha})^{b}\left(\frac{c}{y-\text{alpha}\,y}\right)^{c}\left(-4+c+2c^{2}+c^{3}+\text{alpha}\,(2+c)\,y\right)}{(-1+c)\,(1+c)\,(2+c)}$$

```
(* end outer second integrals *)

(* Step 3: density described in equation  (6) *)
dens1alpha =
 Simplify[(y - b) / y * Exp[ - (y - b - b * Log[(y / b)])]]
(* density if $f_1(0)=b \leq y \leq b/(1-alpha) f_1(alpha)$;
this is relevant for the cases 1 and 2*)
dens2alpha = Simplify[(alpha * Exp[
     - ((b / (1 - alpha) - b - b * Log[(1 / (1 - alpha))]) + ((y - b / (1 - alpha)) alpha))])]
(* density if $y > b/(1-alpha)$, this is relevant for the cases 2 and 4 *)
```

$$\frac{e^{b-y}\left(\frac{y}{b}\right)^{b}(-b+y)}{y}$$

$$\left(\frac{1}{1-\text{alpha}}\right)^{b}\text{alpha}\,e^{-\text{alpha}\,y}$$

```
(* check wheter integral over density is 1 *)
```

```mathematica
Fdens1 = Simplify[Integrate[dens1alpha, {y, b, b / (1 - alpha)},
     Assumptions → {0 < alpha < t < 1, b > 0, c > 0, y > b}]];
Fdens2 = Simplify[Integrate[dens2alpha, {y, b / (1 - alpha), ∞},
     Assumptions → {0 < alpha < t < 1, b > 0, c > 0, y > b}]];
Checkdens = Fdens1 + Fdens2
```

1

```mathematica
(* end density *)

(* Step 4: Obtain components of the value function *)

(* Value function components for $t \leq \alpha$;
since these terms do not depend on $f_2$ we can consider the cases 1,
3 as well as 2,4 jointly *)

V1case13 =
 Simplify[Integrate[dens1alpha * (Outer1integralcaseI1a + Outer1integralcaseI1b),
   {y, b, b / (1 - alpha)}, Assumptions → {0 < alpha < t < 1, b < y < b / (1 - alpha), b > 0}]]
```

$$\frac{1}{2}\left(1 + b - \frac{e^{-\text{alpha}\,b}}{\text{alpha}^2} + e^{\frac{\text{alpha}\,b}{-1+\text{alpha}}} - 2\,(1-\text{alpha})^{-b}\,e^{\frac{\text{alpha}\,b}{-1+\text{alpha}}} + \frac{e^{\frac{\text{alpha}\,b}{-1+\text{alpha}}}}{\text{alpha}^2} + b\,e^{\frac{\text{alpha}\,b}{-1+\text{alpha}}} - \right.$$
$$(1-\text{alpha})^{-b}\,b\,e^{\frac{\text{alpha}\,b}{-1+\text{alpha}}} - b\,e^b\,\text{ExpIntegralEi}[-b] + b\,e^b\,\text{ExpIntegralEi}\left[\frac{b}{-1+\text{alpha}}\right] -$$
$$\left.b\,\text{ExpIntegralEi}[-\text{alpha}\,b] + b\,\text{ExpIntegralEi}\left[\frac{\text{alpha}\,b}{-1+\text{alpha}}\right]\right)$$

```mathematica
V1case24 =
 Simplify[Integrate[dens2alpha * (Outer1integralcaseI24), {y, b / (1 - alpha), ∞},
   Assumptions → {0 < alpha < t < 1, y > 0, b / (1 - alpha) < y, b > 0}]]
```

$(1-\text{alpha})^{-b}\,e^{\frac{\text{alpha}\,b}{-1+\text{alpha}}}\,\text{Outer1integralcaseI24}$

```mathematica
V1case24test =
 Simplify[(Outer1integralcaseI24) * Integrate[dens2alpha, {y, b / (1 - alpha), ∞},
    Assumptions → {0 < alpha < t < 1, y > 0, b / (1 - alpha) < y, b > 0}]]
(* test/simplification, inner integral does not depend on $y$ *)
```

$(1-\text{alpha})^{-b}\,e^{\frac{\text{alpha}\,b}{-1+\text{alpha}}}\,\text{Outer1integralcaseI24}$

```mathematica
(* end value function components for $t \leq \alpha$ *)

(* begin value function components for $t > \alpha$ *)
(* Remark: to integrate over $y$ in the integral associated with $t > \alpha$,
we have to distiguish between the cases where $ b \geq c$ and $b < c$. For the
  first case the abbreviation bgc (g for the german word groesser) is used,
while for the latter case we use bkc (k for the german word kleiner)  *)

V21case1bgc =
 Simplify[Integrate[dens1alpha * (Outer21integralcase1), {y, b, c / (1 - alpha)},
   Assumptions → {0 < alpha < t < 1, y > 0, b < y < b / (1 - alpha), b > c > 0,
     b > 0, b < c / (1 - alpha)}]] (* b \geq c and c/(1-alpha)>b *)
```

$$\frac{1}{2\,\text{alpha}^3\,(1+c)\,(2+c)}\,e^{-\text{alpha}\,b}\left((-1+\text{alpha})\,(-2\,(1+c) + \text{alpha}\,(-2 + (-1+\text{alpha})\,b\,(1+c))) - \right.$$
$$\left(2\,(1+c) + \text{alpha}\,\left(2\,c^2 - b\,(1+c) + \text{alpha}^2\,b\,(3+c\,(2+c)) + \right.\right.$$
$$\left.\left.\left.\text{alpha}\,\left(-2 - b\,\left(2+c+c^2\right) + c\,\left(2+c+c^2\right)\right)\right)\right)\,e^{\text{alpha}\,\left(b+\frac{c}{-1+\text{alpha}}\right)}\right)$$

```
V23case1bgc = Simplify[
  Integrate[dens1alpha * (Outer23integralcase1), {y, b, c / (1 - alpha)}, Assumptions →
    {0 < alpha < t < 1, y > 0, b < y < b / (1 - alpha), b > c > 0, b > 0, b < c / (1 - alpha)}]]
```

$$\frac{1}{2\,(-1+\text{alpha})\,\text{alpha}^3\,(-1+c)\,(1+c)\,(2+c)}$$

$$e^{-\text{alpha}\,b}\left(-2 + 2\,\text{alpha} + 6\,\text{alpha}^2 - 6\,\text{alpha}^3 - \text{alpha}\,b + 5\,\text{alpha}^2\,b - 7\,\text{alpha}^3\,b + 3\,\text{alpha}^4\,b + \right.$$
$$2\,\text{alpha}\,c - 2\,\text{alpha}^3\,c - \text{alpha}^2\,b\,c + 2\,\text{alpha}^3\,b\,c - \text{alpha}^4\,b\,c + 2\,c^2 - 4\,\text{alpha}\,c^2 - \text{alpha}^2\,c^2 +$$
$$\text{alpha}^3\,c^2 + \text{alpha}\,b\,c^2 - 4\,\text{alpha}^2\,b\,c^2 + 5\,\text{alpha}^3\,b\,c^2 - 2\,\text{alpha}^4\,b\,c^2 - 4\,\text{alpha}^2\,c^3 +$$
$$\text{alpha}^3\,c^3 - \text{alpha}^2\,c^4 + 2\,e^{\text{alpha}\,\left(b + \frac{c}{-1+\text{alpha}}\right)} - 2\,\text{alpha}\,e^{\text{alpha}\,\left(b + \frac{c}{-1+\text{alpha}}\right)} - 6\,\text{alpha}^2\,e^{\text{alpha}\,\left(b + \frac{c}{-1+\text{alpha}}\right)} +$$
$$6\,\text{alpha}^3\,e^{\text{alpha}\,\left(b + \frac{c}{-1+\text{alpha}}\right)} - \text{alpha}\,b\,e^{\text{alpha}\,\left(b + \frac{c}{-1+\text{alpha}}\right)} - 3\,\text{alpha}^2\,b\,e^{\text{alpha}\,\left(b + \frac{c}{-1+\text{alpha}}\right)} +$$
$$5\,\text{alpha}^3\,b\,e^{\text{alpha}\,\left(b + \frac{c}{-1+\text{alpha}}\right)} - \text{alpha}^4\,b\,e^{\text{alpha}\,\left(b + \frac{c}{-1+\text{alpha}}\right)} + 4\,\text{alpha}^3\,c\,e^{\text{alpha}\,\left(b + \frac{c}{-1+\text{alpha}}\right)} +$$
$$2\,\text{alpha}^2\,b\,c\,e^{\text{alpha}\,\left(b + \frac{c}{-1+\text{alpha}}\right)} - 4\,\text{alpha}^3\,b\,c\,e^{\text{alpha}\,\left(b + \frac{c}{-1+\text{alpha}}\right)} + 2\,\text{alpha}^4\,b\,c\,e^{\text{alpha}\,\left(b + \frac{c}{-1+\text{alpha}}\right)} -$$
$$2\,c^2\,e^{\text{alpha}\,\left(b + \frac{c}{-1+\text{alpha}}\right)} + 4\,\text{alpha}\,c^2\,e^{\text{alpha}\,\left(b + \frac{c}{-1+\text{alpha}}\right)} + \text{alpha}\,b\,c^2\,e^{\text{alpha}\,\left(b + \frac{c}{-1+\text{alpha}}\right)} -$$
$$\text{alpha}^4\,b\,c^2\,e^{\text{alpha}\,\left(b + \frac{c}{-1+\text{alpha}}\right)} - 2\,\text{alpha}\,c^3\,e^{\text{alpha}\,\left(b + \frac{c}{-1+\text{alpha}}\right)} + 6\,\text{alpha}^2\,c^3\,e^{\text{alpha}\,\left(b + \frac{c}{-1+\text{alpha}}\right)} -$$
$$4\,\text{alpha}^3\,c^3\,e^{\text{alpha}\,\left(b + \frac{c}{-1+\text{alpha}}\right)} + \text{alpha}^2\,b\,c^3\,e^{\text{alpha}\,\left(b + \frac{c}{-1+\text{alpha}}\right)} - \text{alpha}^3\,b\,c^3\,e^{\text{alpha}\,\left(b + \frac{c}{-1+\text{alpha}}\right)} +$$
$$\text{alpha}^3\,b\,\left(2 + 3\,c + c^2\right)\left(4 - 3\,c - c^2 + \text{alpha}\,(-4 + 3\,c)\right)e^{\text{alpha}\,b}\,\text{ExpIntegralEi}[-\text{alpha}\,b] -$$
$$\left. \text{alpha}^3\,b\,\left(4 - 4\,\text{alpha} - 3\,c + 3\,\text{alpha}\,c - c^2\right)\left(2 + 3\,c + c^2\right)e^{\text{alpha}\,b}\,\text{ExpIntegralEi}\left[\frac{\text{alpha}\,c}{-1+\text{alpha}}\right]\right)$$

```
V21case1bkc =
 Simplify[Integrate[dens1alpha * (Outer21integralcase1), {y, b, b / (1 - alpha)},
   Assumptions → {0 < alpha < t < 1, y > b, b < y < b / (1 - alpha), b < c, b > 0, c > 0}]]
```

$$\left(e^{-\text{alpha}\,b}\left((-1+\text{alpha})(-2(1+c) + \text{alpha}(-2 + (-1+\text{alpha})\,b\,(1+c))) - \right.\right.$$
$$\left(2(1+c) + \text{alpha}\left(-2\,\text{alpha} + b + (-2+b)\,c + \text{alpha}^2\,b\,(3 + b + c + b\,c)\right)\right)$$
$$\left.\left.e^{\frac{\text{alpha}^2\,b}{-1+\text{alpha}}}\right)\right) \Big/ \left(2\,\text{alpha}^3\,(1+c)(2+c)\right)$$

```
V21case1bkcacasecsmall =
 Simplify[Integrate[dens1alpha * (Outer21integralcase1csmall), {y, b, b / (1 - alpha)},
   Assumptions → {0 < alpha < t < 1, y > b, b < y < b / (1 - alpha), b < c, b > 0, c > 0}]]
```

$$\frac{1}{2\,\text{alpha}\,(1+c)\,(2+c)}$$

$$(1-\text{alpha})^{-2c}\,e^{-\text{alpha}\,b}\left(4\left((-1+\text{alpha})^2\right)^c + 8\left((-1+\text{alpha})^2\right)^c c + 5\left((-1+\text{alpha})^2\right)^c c^2 + \right.$$

$$\left((-1+\text{alpha})^2\right)^c c^3 - 4\,(1-\text{alpha})^{2c}\,e^{\frac{\text{alpha}^2\,b}{-1+\text{alpha}}} - 8\,(1-\text{alpha})^{2c}\,c\,e^{\frac{\text{alpha}^2\,b}{-1+\text{alpha}}} - $$

$$5\,(1-\text{alpha})^{2c}\,c^2\,e^{\frac{\text{alpha}^2\,b}{-1+\text{alpha}}} - (1-\text{alpha})^{2c}\,c^3\,e^{\frac{\text{alpha}^2\,b}{-1+\text{alpha}}} + $$

$$\left((-1+\text{alpha})^2\right)^c \text{alpha}\,b\,(1+c)\,(2+c)^2\,e^{\text{alpha}\,b}\,\text{ExpIntegralEi}[-\text{alpha}\,b] - $$

$$(1-\text{alpha})^{2c}\,\text{alpha}\,b\,(1+c)\,(2+c)^2\,e^{\text{alpha}\,b}\,\text{ExpIntegralEi}\left[\frac{\text{alpha}\,b}{-1+\text{alpha}}\right] - $$

$$4\,(-(-1+\text{alpha})\,\text{alpha}\,c)^c\,e^{\text{alpha}\,b}\,\text{Gamma}[1-c,\,\text{alpha}\,b] - $$

$$4\,c\,(-(-1+\text{alpha})\,\text{alpha}\,c)^c\,e^{\text{alpha}\,b}\,\text{Gamma}[1-c,\,\text{alpha}\,b] - $$

$$2\,c^2\,(-(-1+\text{alpha})\,\text{alpha}\,c)^c\,e^{\text{alpha}\,b}\,\text{Gamma}[1-c,\,\text{alpha}\,b] + $$

$$4\,(-(-1+\text{alpha})\,\text{alpha}\,c)^c\,e^{\text{alpha}\,b}\,\text{Gamma}\left[1-c,\,\frac{\text{alpha}\,b}{1-\text{alpha}}\right] + $$

$$4\,c\,(-(-1+\text{alpha})\,\text{alpha}\,c)^c\,e^{\text{alpha}\,b}\,\text{Gamma}\left[1-c,\,\frac{\text{alpha}\,b}{1-\text{alpha}}\right] + $$

$$2\,c^2\,(-(-1+\text{alpha})\,\text{alpha}\,c)^c\,e^{\text{alpha}\,b}\,\text{Gamma}\left[1-c,\,\frac{\text{alpha}\,b}{1-\text{alpha}}\right] + $$

$$4\,\text{alpha}\,b\,(-(-1+\text{alpha})\,\text{alpha}\,c)^c\,e^{\text{alpha}\,b}\,\text{Gamma}[-c,\,\text{alpha}\,b] + $$

$$4\,\text{alpha}\,b\,c\,(-(-1+\text{alpha})\,\text{alpha}\,c)^c\,e^{\text{alpha}\,b}\,\text{Gamma}[-c,\,\text{alpha}\,b] + $$

$$2\,\text{alpha}\,b\,c^2\,(-(-1+\text{alpha})\,\text{alpha}\,c)^c\,e^{\text{alpha}\,b}\,\text{Gamma}[-c,\,\text{alpha}\,b] - $$

$$4\,\text{alpha}\,b\,(-(-1+\text{alpha})\,\text{alpha}\,c)^c\,e^{\text{alpha}\,b}\,\text{Gamma}\left[-c,\,\frac{\text{alpha}\,b}{1-\text{alpha}}\right] - $$

$$4\,\text{alpha}\,b\,c\,(-(-1+\text{alpha})\,\text{alpha}\,c)^c\,e^{\text{alpha}\,b}\,\text{Gamma}\left[-c,\,\frac{\text{alpha}\,b}{1-\text{alpha}}\right] - $$

$$\left.2\,\text{alpha}\,b\,c^2\,(-(-1+\text{alpha})\,\text{alpha}\,c)^c\,e^{\text{alpha}\,b}\,\text{Gamma}\left[-c,\,\frac{\text{alpha}\,b}{1-\text{alpha}}\right]\right)$$

```
V23case1bkca =
 Simplify[Integrate[dens1alpha * (Outer23integralcase1), {y, b, b / (1 - alpha)},
    Assumptions → {0 < alpha < t < 1, y > b, b < y < b / (1 - alpha), b < c, b > 0, c > 0}]]
```

$$\frac{1}{2(-1+\text{alpha})\,\text{alpha}^3(-1+c)(1+c)(2+c)}$$

$$e^{-\text{alpha}\,b}\left(e^{\frac{\text{alpha}^2 b}{-1+\text{alpha}}}\left(\text{alpha}\,(-1+c)\,(-2+b-4\,c+b\,c)\left(-1+e^{\frac{\text{alpha}^2 b}{1-\text{alpha}}}\right)+2\left(-1+c^2\right)\left(-1+e^{\frac{\text{alpha}^2 b}{1-\text{alpha}}}\right)+\right.\right.$$

$$\text{alpha}^4\,b\,(-1+c)\left(1-2\,c^2+b\,(3+2\,c)-3\,e^{\frac{\text{alpha}^2 b}{1-\text{alpha}}}-2\,c\left(2+e^{\frac{\text{alpha}^2 b}{1-\text{alpha}}}\right)\right)+\text{alpha}^3\left(-b^2\left(-1+c^2\right)+\right.$$

$$\left.\left(-6-2\,c+c^2+c^3\right)\left(-1+e^{\frac{\text{alpha}^2 b}{1-\text{alpha}}}\right)+b\,(-1+c)\left(-7-3\,c+7\,e^{\frac{\text{alpha}^2 b}{1-\text{alpha}}}+5\,c\,e^{\frac{\text{alpha}^2 b}{1-\text{alpha}}}\right)\right)-$$

$$\left.\text{alpha}^2\,(-1+c)\left(\left(6+6\,c+5\,c^2+c^3\right)\left(-1+e^{\frac{\text{alpha}^2 b}{1-\text{alpha}}}\right)+b\left(-3+5\,e^{\frac{\text{alpha}^2 b}{1-\text{alpha}}}+c\left(-2+4\,e^{\frac{\text{alpha}^2 b}{1-\text{alpha}}}\right)\right)\right)\right)+$$

$$\text{alpha}^3\,b\,\left(2+3\,c+c^2\right)\left(4-3\,c-c^2+\text{alpha}\,(-4+3\,c)\right)\,e^{\text{alpha}\,b}\,\text{ExpIntegralEi}[-\text{alpha}\,b]-$$

$$\left.\text{alpha}^3\,b\,\left(2+3\,c+c^2\right)\left(4-3\,c-c^2+\text{alpha}\,(-4+3\,c)\right)\,e^{\text{alpha}\,b}\,\text{ExpIntegralEi}\left[\frac{\text{alpha}\,b}{-1+\text{alpha}}\right]\right)$$

```
V23case1bkcacasecsmall =
 Simplify[Integrate[dens1alpha * (Outer23integralcase1csmall), {y, b, b / (1 - alpha)},
    Assumptions → {0 < alpha < t < 1, y > b, b < y < b / (1 - alpha), b < c, b > 0, c > 0}]]
```

$$\frac{1}{\text{alpha}\,(-1+c)(1+c)(2+c)}$$

$$\left(-\frac{\text{alpha}\,c}{-1+\text{alpha}}\right)^c\left(\left(-4+c+2\,c^2+c^3-\text{alpha}\,b\,(2+c)\right)\text{Gamma}[1-c,\text{alpha}\,b]+\right.$$

$$\left(4-c-2\,c^2-c^3+\text{alpha}\,b\,(2+c)\right)\text{Gamma}\left[1-c,\frac{\text{alpha}\,b}{1-\text{alpha}}\right]+2\,\text{Gamma}[2-c,\text{alpha}\,b]+$$

$$c\,\text{Gamma}[2-c,\text{alpha}\,b]-2\,\text{Gamma}\left[2-c,\frac{\text{alpha}\,b}{1-\text{alpha}}\right]-c\,\text{Gamma}\left[2-c,\frac{\text{alpha}\,b}{1-\text{alpha}}\right]+$$

$$4\,\text{alpha}\,b\,\text{Gamma}[-c,\text{alpha}\,b]-\text{alpha}\,b\,c\,\text{Gamma}[-c,\text{alpha}\,b]-$$

$$2\,\text{alpha}\,b\,c^2\,\text{Gamma}[-c,\text{alpha}\,b]-\text{alpha}\,b\,c^3\,\text{Gamma}[-c,\text{alpha}\,b]-$$

$$4\,\text{alpha}\,b\,\text{Gamma}\left[-c,\frac{\text{alpha}\,b}{1-\text{alpha}}\right]+\text{alpha}\,b\,c\,\text{Gamma}\left[-c,\frac{\text{alpha}\,b}{1-\text{alpha}}\right]+$$

$$\left.2\,\text{alpha}\,b\,c^2\,\text{Gamma}\left[-c,\frac{\text{alpha}\,b}{1-\text{alpha}}\right]+\text{alpha}\,b\,c^3\,\text{Gamma}\left[-c,\frac{\text{alpha}\,b}{1-\text{alpha}}\right]\right)$$

```
V21case2bkca1 = Simplify[Integrate[dens2alpha * (Outer21integralcase2),
    {y, b / (1 - alpha), 1}, Assumptions → {0 < alpha < 1, c > b > 0}]]
```

$$\frac{1}{2(1+c)(2+c)}$$

$$\left(-1+\frac{1}{\text{alpha}}\right)^{-c} e^{\frac{\text{alpha}\,(1-\text{alpha}+b)}{-1+\text{alpha}}}\left(\left(-1+\frac{1}{\text{alpha}}\right)^c(1+c)(2+c)^2\left(e^{\text{alpha}}-e^{\frac{\text{alpha}\,b}{1-\text{alpha}}}\right)+2\,c^c\left(2+2\,c+c^2\right)\right.$$

$$\left.e^{\text{alpha}+\frac{\text{alpha}\,b}{1-\text{alpha}}}\text{Gamma}[1-c,\text{alpha}]-2\,c^c\left(2+2\,c+c^2\right)e^{\text{alpha}+\frac{\text{alpha}\,b}{1-\text{alpha}}}\text{Gamma}\left[1-c,\frac{\text{alpha}\,b}{1-\text{alpha}}\right]\right)$$

```
V21case2bkca2 = Simplify[Integrate[dens2alpha * (Outer21integralcase2),
   {y, c / (1 - alpha), 1}, Assumptions → {0 < alpha < 1, c > b > 0}]]
```

$$\frac{1}{2(1+c)(2+c)}\left(-1+\frac{1}{alpha}\right)^{-c} e^{\frac{alpha(1-alpha+c)}{-1+alpha}}$$

$$\left(\left(-1+\frac{1}{alpha}\right)^c (1+c)(2+c)^2 \left(e^{alpha} - e^{\frac{alpha\, c}{1-alpha}}\right) + 2\, c^c \left(2+2c+c^2\right) e^{alpha-\frac{alpha\, c}{-1+alpha}}\right.$$

$$\left. \text{Gamma}[1-c,\, alpha] - 2\, c^c \left(2+2c+c^2\right) e^{alpha-\frac{alpha\, c}{-1+alpha}} \text{Gamma}\!\left[1-c,\, -\frac{alpha\, c}{-1+alpha}\right]\right)$$

```
V21case2bkc = Simplify[V21case2bkca1 - V21case2bkca2]
```

$$\frac{1}{2(1+c)(2+c)}$$

$$\left(-1+\frac{1}{alpha}\right)^{-c}\left(e^{\frac{alpha(1-alpha+b)}{-1+alpha}}\left(\left(-1+\frac{1}{alpha}\right)^c (1+c)(2+c)^2 \left(e^{alpha}-e^{\frac{alpha\, b}{1-alpha}}\right) + 2\, c^c \left(2+2c+c^2\right)\right.\right.$$

$$\left.e^{alpha+\frac{alpha\, b}{1-alpha}} \text{Gamma}[1-c,\, alpha] - 2\, c^c \left(2+2c+c^2\right) e^{alpha+\frac{alpha\, b}{1-alpha}} \text{Gamma}\!\left[1-c,\, \frac{alpha\, b}{1-alpha}\right]\right) -$$

$$e^{\frac{alpha(1-alpha+c)}{-1+alpha}}\left(\left(-1+\frac{1}{alpha}\right)^c (1+c)(2+c)^2 \left(e^{alpha}-e^{\frac{alpha\, c}{1-alpha}}\right) + 2\, c^c \left(2+2c+c^2\right) e^{alpha-\frac{alpha\, c}{-1+alpha}}\right.$$

$$\left.\left. \text{Gamma}[1-c,\, alpha] - 2\, c^c \left(2+2c+c^2\right) e^{alpha-\frac{alpha\, c}{-1+alpha}} \text{Gamma}\!\left[1-c,\, -\frac{alpha\, c}{-1+alpha}\right]\right)\right)$$

```
V23case2bkc = Simplify[
  Integrate[dens2alpha * (Outer23integralcase2), {y, b / (1 - alpha), c / (1 - alpha)},
   Assumptions → {0 < alpha < t < 1, y > b, b / (1 - alpha) < y, b < c, b > 0, c > 0}]]
```

$$\frac{1}{2(-1+alpha)\, alpha^2\, (-1+c)(1+c)(2+c)} e^{\frac{alpha(b+c)}{-1+alpha}}$$

$$\left(-2\left(-1+c^2+alpha\left(1+(-2+c)c^2\right) - 3\, alpha^2\left(-1+c^3\right) + alpha^3\left(-3-2c+2c^3\right)\right) e^{\frac{alpha\, b}{1-alpha}} - \right.$$

$$\left(2 - 2c^2 - 2\, alpha\, (-1+c)(-1+b+(-2+b)c) + \right.$$

$$\quad alpha^2 \left(-6+b^2-2bc+(1-(-2+b)b)c^2+4c^3+c^4\right) +$$

$$\quad \left.\left. alpha^3 \left(6-c\left(-2+c+c^2\right)+b^2\left(-3+c+2c^2\right)-2b(-1+c)(1+c(3+c))\right)\right) e^{-\frac{alpha\, c}{-1+alpha}}\right)$$

```
V21case3bgc = Simplify[Integrate[dens1alpha * (Outer21integralcase3),
    {y, c / (1 - alpha), b / (1 - alpha)}, Assumptions → {0 < alpha < 1, y > b > c > 0}]]
```

$$\frac{1}{2\,\text{alpha}\,(1+c)\,(2+c)}$$

$$\left(-4\,e^{\frac{\text{alpha}\,b}{-1+\text{alpha}}} - 8\,c\,e^{\frac{\text{alpha}\,b}{-1+\text{alpha}}} - 5\,c^2\,e^{\frac{\text{alpha}\,b}{-1+\text{alpha}}} - c^3\,e^{\frac{\text{alpha}\,b}{-1+\text{alpha}}} + 4\,e^{\frac{\text{alpha}\,c}{-1+\text{alpha}}} + 8\,c\,e^{\frac{\text{alpha}\,c}{-1+\text{alpha}}} + 5\,c^2\,e^{\frac{\text{alpha}\,c}{-1+\text{alpha}}} +\right.$$

$$c^3\,e^{\frac{\text{alpha}\,c}{-1+\text{alpha}}} - \text{alpha}\,b\,(1+c)\,(2+c)^2\,\text{ExpIntegralEi}\left[\frac{\text{alpha}\,b}{-1+\text{alpha}}\right] +$$

$$\text{alpha}\,b\,(1+c)\,(2+c)^2\,\text{ExpIntegralEi}\left[\frac{\text{alpha}\,c}{-1+\text{alpha}}\right] +$$

$$4\left(\frac{\text{alpha}\,c}{1-\text{alpha}}\right)^c \text{Gamma}\left[1-c,\,\frac{\text{alpha}\,b}{1-\text{alpha}}\right] + 4\,c\left(\frac{\text{alpha}\,c}{1-\text{alpha}}\right)^c \text{Gamma}\left[1-c,\,\frac{\text{alpha}\,b}{1-\text{alpha}}\right] +$$

$$2\,c^2\left(\frac{\text{alpha}\,c}{1-\text{alpha}}\right)^c \text{Gamma}\left[1-c,\,\frac{\text{alpha}\,b}{1-\text{alpha}}\right] - 4\left(\frac{\text{alpha}\,c}{1-\text{alpha}}\right)^c \text{Gamma}\left[1-c,\,\frac{\text{alpha}\,c}{1-\text{alpha}}\right] -$$

$$4\,c\left(\frac{\text{alpha}\,c}{1-\text{alpha}}\right)^c \text{Gamma}\left[1-c,\,\frac{\text{alpha}\,c}{1-\text{alpha}}\right] - 2\,c^2\left(\frac{\text{alpha}\,c}{1-\text{alpha}}\right)^c \text{Gamma}\left[1-c,\,\frac{\text{alpha}\,c}{1-\text{alpha}}\right] -$$

$$4\,\text{alpha}\,b\left(\frac{\text{alpha}\,c}{1-\text{alpha}}\right)^c \text{Gamma}\left[-c,\,\frac{\text{alpha}\,b}{1-\text{alpha}}\right] - 4\,\text{alpha}\,b\,c\left(\frac{\text{alpha}\,c}{1-\text{alpha}}\right)^c \text{Gamma}\left[-c,\,\frac{\text{alpha}\,b}{1-\text{alpha}}\right] -$$

$$2\,\text{alpha}\,b\,c^2\left(\frac{\text{alpha}\,c}{1-\text{alpha}}\right)^c \text{Gamma}\left[-c,\,\frac{\text{alpha}\,b}{1-\text{alpha}}\right] +$$

$$4\,\text{alpha}\,b\left(\frac{\text{alpha}\,c}{1-\text{alpha}}\right)^c \text{Gamma}\left[-c,\,\frac{\text{alpha}\,c}{1-\text{alpha}}\right] + 4\,\text{alpha}\,b\,c\left(\frac{\text{alpha}\,c}{1-\text{alpha}}\right)^c \text{Gamma}\left[-c,\,\frac{\text{alpha}\,c}{1-\text{alpha}}\right] +$$

$$\left.2\,\text{alpha}\,b\,c^2\left(\frac{\text{alpha}\,c}{1-\text{alpha}}\right)^c \text{Gamma}\left[-c,\,\frac{\text{alpha}\,c}{1-\text{alpha}}\right]\right)$$

```
V23case3bgc = Simplify[
    Integrate[dens1alpha * (Outer23integralcase3), {y, c / (1 - alpha), b / (1 - alpha)},
     Assumptions → {0 < alpha < 1, y > b, b > c, b > 0, c > 0}]]
```

$$\frac{1}{\text{alpha}\,(-1+c)\,(1+c)\,(2+c)}$$

$$\left(\frac{\text{alpha}\,c}{1-\text{alpha}}\right)^c \left(\left(4-c-2\,c^2-c^3+\text{alpha}\,b\,(2+c)\right)\text{Gamma}\left[1-c,\,\frac{\text{alpha}\,b}{1-\text{alpha}}\right] +\right.$$

$$\left(-4+c+2\,c^2+c^3-\text{alpha}\,b\,(2+c)\right)\text{Gamma}\left[1-c,\,\frac{\text{alpha}\,c}{1-\text{alpha}}\right] - 2\,\text{Gamma}\left[2-c,\,\frac{\text{alpha}\,b}{1-\text{alpha}}\right] -$$

$$c\,\text{Gamma}\left[2-c,\,\frac{\text{alpha}\,b}{1-\text{alpha}}\right] + 2\,\text{Gamma}\left[2-c,\,\frac{\text{alpha}\,c}{1-\text{alpha}}\right] + c\,\text{Gamma}\left[2-c,\,\frac{\text{alpha}\,c}{1-\text{alpha}}\right] -$$

$$4\,\text{alpha}\,b\,\text{Gamma}\left[-c,\,\frac{\text{alpha}\,b}{1-\text{alpha}}\right] + \text{alpha}\,b\,c\,\text{Gamma}\left[-c,\,\frac{\text{alpha}\,b}{1-\text{alpha}}\right] +$$

$$2\,\text{alpha}\,b\,c^2\,\text{Gamma}\left[-c,\,\frac{\text{alpha}\,b}{1-\text{alpha}}\right] + \text{alpha}\,b\,c^3\,\text{Gamma}\left[-c,\,\frac{\text{alpha}\,b}{1-\text{alpha}}\right] +$$

$$4\,\text{alpha}\,b\,\text{Gamma}\left[-c,\,\frac{\text{alpha}\,c}{1-\text{alpha}}\right] - \text{alpha}\,b\,c\,\text{Gamma}\left[-c,\,\frac{\text{alpha}\,c}{1-\text{alpha}}\right] -$$

$$\left.2\,\text{alpha}\,b\,c^2\,\text{Gamma}\left[-c,\,\frac{\text{alpha}\,c}{1-\text{alpha}}\right] - \text{alpha}\,b\,c^3\,\text{Gamma}\left[-c,\,\frac{\text{alpha}\,c}{1-\text{alpha}}\right]\right)$$

**V21case4bgc = Simplify[Integrate[dens2alpha \* (Outer21integralcase4),**
  **{y, b / (1 - alpha), ∞}, Assumptions → {0 < alpha < 1, y > b, b > c > 0}]]**

$$\frac{1}{2} (2+c) e^{\frac{alpha\, b}{-1+alpha}} + \frac{alpha\, b \left(\frac{c}{b}\right)^c (2+c\,(2+c))\, \text{ExpIntegralE}\left[c, \frac{alpha\, b}{1-alpha}\right]}{(-1+alpha)\,(1+c)\,(2+c)}$$

**V23case4bgc = Simplify[Integrate[dens2alpha \* (Outer23integralcase4),**
  **{y, b / (1 - alpha), ∞}, Assumptions → {0 < alpha < t < 1, y > 0, b > c, b > 1, c > 1}]]**

$$\frac{1}{(-1+c)(1+c)(2+c)}$$

$$\text{alpha}\left(\frac{\text{alpha}\, b^2 \left(\frac{c}{b}\right)^c (2+c)\, \Gamma[-2+c]}{(-1+alpha)^2\, \Gamma[-1+c]} - \frac{b\left(\frac{c}{b}\right)^c (-4+c+2c^2+c^3)\, \Gamma[-1+c]}{(-1+alpha)\, \Gamma[c]} + \frac{1}{alpha}\right.$$

$$(1-alpha)^{-1-2c}\,(alpha\, b\, c)^c\, (4+3c+c^2)\left(alpha\,\left(\frac{(-1+alpha)^2}{alpha\, b^2}\right)^c b\left(-1+e^{\frac{alpha\, b}{-1+alpha}}\right)+\right.$$

$$(-1+alpha)\left(\frac{1-alpha}{b}\right)^c \Gamma\left[2-c, \frac{alpha\, b}{1-alpha}\right]\right)-\frac{1}{alpha}$$

$$\left(\frac{c}{1-alpha}\right)^c (2+c)\left(-alpha^c\, \Gamma[2-c]+\left(b^{-c}\left(-(1-alpha)^c\, alpha^2\, b^2\left(-1+e^{\frac{alpha\, b}{-1+alpha}}\right)+\right.\right.\right.$$

$$(-1+alpha)^2\,(alpha\, b)^c\,(-2+c)\,\Gamma[2-c]+$$

$$\left.\left.\left.\left.(-1+alpha)^2\,(alpha\, b)^c\, \Gamma\left[3-c, \frac{alpha\, b}{1-alpha}\right]\right)\right)\middle/\left((-1+alpha)^2\,(-2+c)\right)\right)\right)$$

**V21case4bkc = Simplify[Integrate[dens2alpha \* Outer21integralcase4,**
  **{y, c / (1 - alpha), ∞}, Assumptions → {0 < alpha < 1, y > c / (1 - alpha), c > b > 0}]]**

$$\frac{1}{2}(2+c)\,e^{\frac{alpha\, c}{-1+alpha}} + \frac{alpha\, c\,(2+c\,(2+c))\,\text{ExpIntegralE}\left[c, \frac{alpha\, c}{1-alpha}\right]}{(-1+alpha)(1+c)(2+c)}$$

```
V23case4bkc = Simplify[Integrate[dens2alpha * (Outer23integralcase4),
   {y, c / (1 - alpha), ∞}, Assumptions → {0 < alpha < t < 1, y > 0, c > b, b > 1, c > 1}]]
```

$$\frac{1}{(-1+c)(1+c)(2+c)}$$

$$\text{alpha}\left(\frac{\text{alpha}\, c^2 (2+c)\, \text{Gamma}[-2+c]}{(-1+\text{alpha})^2\, \text{Gamma}[-1+c]} - \frac{c\left(-4+c+2c^2+c^3\right)\text{Gamma}[-1+c]}{(-1+\text{alpha})\, \text{Gamma}[c]} - \right.$$

$$\left.\left(\left(\frac{c}{1-\text{alpha}}\right)^c (4+3c+c^2)\left(\text{alpha}\left(\frac{1-\text{alpha}}{c}\right)^c c\left(-1+e^{\frac{\text{alpha}\, c}{-1+\text{alpha}}}\right) + \right.\right.\right.$$

$$\left.\left.\left.(-1+\text{alpha})\, \text{alpha}^c\, \text{Gamma}\left[2-c, \frac{\text{alpha}\, c}{1-\text{alpha}}\right]\right)\right)\Big/((-1+\text{alpha})\,\text{alpha})\right) -$$

$$\text{alpha}\left(\frac{c}{1-\text{alpha}}\right)^c (2+c)\left(-\text{alpha}^{-2+c}\,\text{Gamma}[2-c] + 1/(-3+c)\,\text{alpha}\left(\frac{1-\text{alpha}}{c}\right)^{-3+c}\right.$$

$$\left(\frac{1-\text{alpha}}{\text{alpha}(-2+c)} + \frac{3(-1+\text{alpha})}{\text{alpha}(-2+c)\, c} + \frac{(-1+\text{alpha})\, e^{\frac{\text{alpha}\, c}{-1+\text{alpha}}}}{\text{alpha}(-2+c)} - \frac{3(-1+\text{alpha})\, e^{\frac{\text{alpha}\, c}{-1+\text{alpha}}}}{\text{alpha}(-2+c)\, c} - \right.$$

$$3\left(\frac{1-\text{alpha}}{\text{alpha}\, c}\right)^{3-c}\text{Gamma}[2-c] + \left(-1+\frac{1}{\text{alpha}}\right)^{3-c} c^{-2+c}\,\text{Gamma}[2-c] -$$

$$\left.\left.\frac{3\left(\frac{1-\text{alpha}}{\text{alpha}\, c}\right)^{3-c}\text{Gamma}\left[3-c, \frac{\text{alpha}\, c}{1-\text{alpha}}\right]}{-2+c} + \frac{\left(-1+\frac{1}{\text{alpha}}\right)^{3-c} c^{-2+c}\,\text{Gamma}\left[3-c, \frac{\text{alpha}\, c}{1-\text{alpha}}\right]}{-2+c}\right)\right)$$

```
V22case24 = Simplify[Integrate[dens2alpha * (Outer22integralcase24),
   {y, b / (1 - alpha), ∞}, Assumptions → {0 < alpha < t < 1, y > 0, c > b, b > 1, c > 1}]]
```

$$\frac{e^{\frac{\text{alpha}\, b}{-1+\text{alpha}}}}{-1+c^2}$$

```
V22case13 = Simplify[Integrate[dens1alpha * (Outer22integralcase13),
   {y, b, b / (1 - alpha)}, Assumptions → {0 < alpha < 1}]]
```

$$\frac{1}{-1+c^2}\left(-\frac{e^{-\text{alpha}\, b}\left(-1+e^{\frac{\text{alpha}^2\, b}{-1+\text{alpha}}}\right)}{\text{alpha}} + b\left(\text{ExpIntegralEi}[-\text{alpha}\, b] - \text{ExpIntegralEi}\left[\frac{\text{alpha}\, b}{-1+\text{alpha}}\right]\right)\right)$$

```
(* end value function, t > \alp1ha *)



(* end value function *)
```

```
(*Some Remarks: In the above value function components we consider the cases 1-
    4 and distiguished between $b \geq c$ and $b > c$. For $b \geq c$,
Case 2 is impossible, while for $b<c$  our Case 3 becomes impossible. By
   this fact to obtain the value function the expressions for the Cases 1,
3,4 are used whenever $b \geq c$, while for $b< c$ the Cases 1,
2,4 are used to obtain the value function. If $b=
 c$ all expressions containing "Case2" or "Case3" are zero. This follows
   dicrectly when looking at the integrals in equation (8). Numerically,
we observe that, using the expressions obtained above,
with $b=c$ these components are zero.    *)
```



\end{document}